\documentclass[a4paper]{article}

\usepackage[latin1]{inputenc}
\usepackage{graphicx}
\usepackage{amssymb,amsfonts,amsmath}
\usepackage{theorem}
\usepackage{hyperref}
\usepackage{graphicx}

\newtheorem{lemma}{Lemma}[section]
\newtheorem{proposition}{Proposition}[section]
{\theorembodyfont{\rmfamily}
\newtheorem{definition}{Definition}[section]
\newtheorem{remark}{Remark}[section]
}

\newenvironment{pol1}[1][Proof of Lemma \ref{lema1}]{\noindent\textbf{#1.} }{\newline \hspace*{\textwidth}\hspace*{-0,4cm} \rule{0.5em}{0.5em} \vspace{0,2cm}}
\newenvironment{pol2}[1][Proof of Lemma \ref{lema2}]{\noindent\textbf{#1.} }{\newline \hspace*{\textwidth}\hspace*{-0,4cm} \rule{0.5em}{0.5em} \vspace{0,2cm}}
\newenvironment{pol3}[1][Proof of Lemma \ref{lema3}]{\noindent\textbf{#1.} }{\newline \hspace*{\textwidth}\hspace*{-0,4cm} \rule{0.5em}{0.5em} \vspace{0,2cm}}
\newenvironment{pot1}[1][Proof of Proposition \ref{prop1}]{\noindent\textbf{#1.} }{\newline \hspace*{\textwidth}\hspace*{-0,4cm} \rule{0.5em}{0.5em} \vspace{0,2cm}}
\newenvironment{pot2}[1][Proof of Proposition \ref{prop2}]{\noindent\textbf{#1.} }{\newline \hspace*{\textwidth}\hspace*{-0,4cm} \rule{0.5em}{0.5em} \vspace{0,2cm}}
\newenvironment{pot3}[1][Proof of Proposition \ref{prop3}]{\noindent\textbf{#1.} }{\newline \hspace*{\textwidth}\hspace*{-0,4cm} \rule{0.5em}{0.5em} \vspace{0,2cm}}
\newenvironment{pot4}[1][Proof of Proposition \ref{prop4}]{\noindent\textbf{#1.} }{\newline \hspace*{\textwidth}\hspace*{-0,4cm} \rule{0.5em}{0.5em} \vspace{0,2cm}}
\newenvironment{pot5}[1][Proof of Proposition \ref{prop5}]{\noindent\textbf{#1.} }{\newline \hspace*{\textwidth}\hspace*{-0,4cm} \rule{0.5em}{0.5em} \vspace{0,2cm}}
\newenvironment{pot6}[1][Proof of Proposition \ref{prop6}]{\noindent\textbf{#1.} }{\newline \hspace*{\textwidth}\hspace*{-0,4cm} \rule{0.5em}{0.5em} \vspace{0,2cm}}

\begin{document}

\title{\textbf{Existence and uniqueness of nontrivial collocation solutions of
  implicitly linear homogeneous Volterra integral equations}}
\author{R. Ben\'{\i}tez$^{1,*}$, V. J. Bol\'os$^{2,*}$ \\
\\
{\small $^1$ Dpto. Matem\'aticas, Centro Universitario de Plasencia, Universidad de Extremadura.}\\
{\small Avda. Virgen del Puerto 2, 10600 Plasencia, Spain.}\\
{\small e-mail\textup{: \texttt{rbenitez@unex.es}}} \\
\\
{\small $^2$ Dpto. Matem\'aticas para la Econom\'{\i}a y la Empresa, Facultad de Econom\'{\i}a,}\\
{\small Universidad de Valencia. Avda. Tarongers s/n, 46022 Valencia, Spain.}\\
{\small e-mail\textup{: \texttt{vicente.bolos@uv.es}}} \\
}
\footnotetext{$^*$Work supported by project MTM2008-05460, Spain.}
\date{September 2010
}

\maketitle

\begin{abstract}
  We analyze collocation methods for nonlinear homogeneous
  Volterra-Hammerstein integral equations with non-Lipschitz
  nonlinearity. We present different kinds of existence and uniqueness
  of nontrivial collocation solutions and we give conditions for such
  existence and uniqueness in some cases. Finally we illustrate these
  methods with an example of a collocation problem, and we give some
  examples of collocation problems that do not fit in the cases
  studied previously.
\end{abstract}

\section{Introduction}

The aim of this paper is the numerical analysis of the nonlinear
homogeneous Volterra-Hammerstein integral equation (HVHIE)
\begin{equation}
\label{eq0}
y\left( t\right) = \left( \mathcal{H}y\right) \left( t\right)
:=\int _0 ^t K\left( t,s\right) G\left( y\left( s\right) \right) \,
\mathrm{d}s , \qquad t\in I:=\left[ 0,T\right] ,
\end{equation}
by means of collocation methods on spaces of local polynomials.
This equation has multiple applications in analysis and physics,
as for example, the study of viscoelastic materials, the renewal
equation, seismic response, transverse oscillations or flows of
heat (see \cite{Mann51,Rahm07}).

Functions $K$ and $G$ are called \textit{kernel} and
\textit{nonlinearity} respectively, and we assume that the
following \textit{general conditions} are always held, even if they
are not explicitly mentioned.

\begin{itemize}
\item\textbf{Over $K$.} The kernel $K:\mathbb{R}^2\to \left[ 0,+\infty
  \right[ $ is a locally bounded function and its support is in
  $\left\{ \left( t,s\right) \in \mathbb{R}^2\, :\, 0\leq s\leq
    t\right\} $.

  For every $t>0$, the map $s\mapsto K\left( t,s\right) $ is locally
  integrable, and $\int _0 ^t K\left( t,s\right) \mathrm{d}s$ is a
  strictly increasing function.

\item \textbf{Over $G$.} The nonlinearity $G:\left[ 0,+\infty \right[
  \to \left[ 0,+\infty \right[ $ is a continuous, strictly increasing
  function, and $G\left( 0\right) =0$.
\end{itemize}

Note that since $G(0)=0$, the zero function is a solution of equation
(\ref{eq0}), known as \textit{trivial solution}, and so, uniqueness of
solutions is no longer a desired property for equation (\ref{eq0})
because we are obviously interested in nontrivial solutions.
Existence and uniqueness of nontrivial solutions of equation
(\ref{eq0}), as well as their properties, have been deeply studied in
a wide range of cases for $K$ and $G$
\cite{Okra91,Szwa92,Arias99,Ari00,AriBen03,AriBenBol05}, especially in
the case of convolution equations, i.e. $K(t,s) = k(t-s)$.
In general, a necessary and sufficient condition for the existence of a nontrivial solution is
the existence of a nontrivial subsolution; that is, a positive function $u$ such that $u(t)\leq \left( \mathcal{H}u\right) (t)$.
So most of the results on existence of nontrivial solutions are, indeed, characterizations
of the existence of subsolutions. For instance, in \cite{Arias99} it can be found the next result:
under the \textit{general conditions}, equation (\ref{eq0}) has a nontrivial solution if and only if
there is a positive integrable function $f(x)$ such that $\int _0^x \mathcal{K}\left( \mathcal{F}(x)-\mathcal{F}(s)\right) \, \mathrm{d}s\geq G^{-1}(x)$, $x\geq 0$, where $\mathcal{K}(x):=\int _0^x k(s)\, \mathrm{d}s$ and $\mathcal{F}(x):=\int _0^x f(s)\, \mathrm{d}s$.

It is important to note that usually, in the analysis of solutions for
non-homogeneous Volterra integral equations (and their numerical
approximations), most of the existence and uniqueness theorems require
that a Lipschitz condition is held by the nonlinearity (with some
exceptions, for instance \cite{Frisch97}). This is not our case, since
it is well known that if the nonlinearity is Lipschitz-continuous,
then the unique solution of (\ref{eq0}) is the trivial one
\cite{Brun04}. Thus, the case we are going to consider in this paper
is beyond the scope of classical results of numerical analysis of
non-homogeneous Volterra integral equations, in the sense that we need a non-Lipschitz nonlinearity.

Actually, there is a wide range of numerical methods available for
solving integral equations (see \cite{Baker96} for a comprehensive
survey on the subject): iterative methods, wavelet methods
\cite{Vainikko05,Yousefi06,Xiao06,Gao08}, generalized Runge-Kutta
methods \cite{Capo07,Calvo07}, or even Monte Carlo methods
\cite{Farnoosh08}. Collocation methods \cite{Brun04,Conte06} have
proved to be very suitable for a wide range of equations, because of
their accuracy, stability and rapid convergence.
In this work we use collocation methods to solve the nonlinear
HVHIE (\ref{eq0}) written in its implicitly linear form (see below). We
also give conditions for different kinds of existence and uniqueness
of nontrivial collocation solutions for the corresponding collocation
equations.

We organize this paper into four sections. In Section \ref{sec:prel},
we write equation (\ref{eq0}) in its implicitly linear form and we
describe the corresponding collocation equations; moreover, we define
the concept of nontrivial collocation solution.  In Section
\ref{sec:exist}, we present different kinds of existence of nontrivial
collocation solutions and we give conditions for their existence and
uniqueness in some cases, considering convolution and nonconvolution
kernels. In Section \ref{sec:examples}, we illustrate the collocation
methods and their numerical convergence with an example, showing how
the errors change as the collocation points vary. Moreover, we give
some examples of collocation problems that do not fit in the cases
studied in the paper.  Finally, we present the proofs of the main
results in an appendix at the end of the paper, for the sake of
readability.

\section{Preliminary concepts}
\label{sec:prel}

Let us consider the nonlinear homogeneous Volterra-Hammerstein
integral equation (HVHIE) given by (\ref{eq0}).
Taking $z:=G\circ y$, equation (\ref{eq0}) can be written as an
\textit{implicitly linear} homogeneous Volterra integral equation
(HVIE) for $z$:
\begin{equation}
\label{eq2}
z(t) = G\left( \left( \mathcal{V}z\right)
  \left( t\right) \right) =G\left( \int _0 ^t K\left( t,s\right)
  z(s) \, \mathrm{d}s \right) , \qquad t\in I,
\end{equation}
where $\mathcal{V}$ is the linear \textit{Volterra operator}. So, if
$z$ is a solution of (\ref{eq2}), then $y:=\mathcal{V}z$ is a solution
of (\ref{eq0}). It is known (see \cite[p. 143]{Kras84}) that, under
suitable assumptions on the nonlinearity $G$, there is a one-to-one
correspondence between solutions of (\ref{eq0}) and
(\ref{eq2}). Particularly, if $G$ is injective then $y=G^{-1}\circ z$
and hence this correspondence is given, which is granted by the
\textit{general conditions} exposed above.

\subsection{Collocation problems for implicitly linear HVIEs}

First, we are going to introduce the collocation problem associated to equation
(\ref{eq2}), and give the equations for determining a collocation solution, that
we will use for approximating a solution of (\ref{eq2}) or (\ref{eq0}) (see Remark \ref{remnlHVIE} below).

Let $I_h:=\left\{ t_n\, :\, \, 0=t_0<t_1<\ldots <t_N=T\right\} $ be a
mesh (not necessarily uniform) on the interval $I=\left[ 0,T\right] $,
and set $\sigma _n:=\left] t_n,t_{n+1}\right] $ with lengths
$h_n:=t_{n+1}-t_n$ $(n=0,\ldots ,N-1)$. The quantity
$h:=\mathrm{max}\left\{ h_n\, :\, \, 0\leq n\leq N-1\right\} $ is
called the \textit{stepsize}.

Given a set of $m$ \textit{collocation parameters} $\left\{c_i\, :\,
  0\leq c_1<\ldots <c_m\leq 1\right\} $, the \textit{collocation
  points} are given by $t_{n,i} :=t_n+c_ih_n$ $(n=0,\ldots ,N-1)$
$(i=1,\ldots ,m)$, and the set of collocation points is denoted by
$X_h$.

All this defines a \textit{collocation problem} for equation
(\ref{eq2}) (see \cite{Brun92}, \cite[p. 117]{Brun04}), and a
\textit{collocation solution} $z_h$ is given by the
\textit{collocation equation}
\begin{equation}
\label{eqcolloc}
z_h\left( t\right) = G\left( \int _0 ^t K\left( t,s\right) z_h\left( s\right)
\, \mathrm{d}s \right) , \qquad t\in X_h\, ,
\end{equation}
where $z_h$ is in the space of piecewise polynomials of degree less
than $m$ (see \cite[p. 85]{Brun04}). Note
that the identically zero function is always a collocation solution,
since $G\left( 0\right) =0$.

\begin{remark}
\label{remnlHVIE}
From now on, a ``collocation problem'' or a ``collocation solution''
will be always referred to the implicitly linear equation
(\ref{eq2}). So, if we want to obtain an estimation of a solution of
the nonlinear HVHIE (\ref{eq0}), then we have to consider
$y_h:=\mathcal{V}z_h$.
\end{remark}

As it is stated in \cite{Brun04}, a collocation solution $z_h$ is
completely determined by the coefficients $Z_{n,i}:=z_h\left(
t_{n,i}\right) $ $(n=0,\ldots ,N-1)$ $(i=1,\ldots ,m)$, since
$z_h\left( t_n+vh_n\right) =\sum _{j=1}^m L_j\left( v\right) Z_{n,j}$
for all $v\in \left] 0,1\right] $, where $L_j\left( v\right) :=\prod
_{k\neq j}^m \frac{v-c_k}{c_j-c_k}$ $(j=1,\ldots ,m$) are the Lagrange
fundamental polynomials with respect to the collocation parameters.
The values of $Z_{n,i}$ are given by the systems
\begin{equation}
\label{eqznigeneral}
Z_{n,i}=G\left( F_n\left( t_{n,i}\right) +h_n\sum _{j=1}^m
B_n\left( i,j\right) Z_{n,j}\right) ,
\end{equation}
where
\begin{equation}
\label{eqbnijgeneral}
B_n\left( i,j\right) := \int _0^{c_i}K\left( t_{n,i},t_n+sh_n\right)
L_j\left( s\right) \, \mathrm{d}s.
\end{equation}
and
\begin{equation}
\label{eqlag0}
F_n\left( t\right) :=\int _0^{t_n} K\left( t,s\right)
z_h\left( s\right) \, \mathrm{d}s.
\end{equation}
The term $F_n\left( t_{n,i}\right) $ is the \textit{lag term}, and can be expressed in the form
\[
F_n\left( t_{n,i}\right) =\sum _{l=0}^{n-1} h_l \sum _{j=1}^m B_n^l
\left( i,j\right) Z_{l,j} ,
\]
where
\[
B_n^l\left( i,j\right) := \int _0^1 K\left( t_{n,i},t_l+sh_l\right)
L_j\left( s\right) \, \mathrm{d}s
\]
with $n=0,\ldots ,N-1$, $l=0,\ldots ,n-1$, $i=1,\ldots ,m$, $j=1,\ldots ,m$.

\begin{remark}
\label{notabnconv}
For convolution kernels, $K\left( t,s\right) =k\left( t-s\right) $,
expression (\ref{eqbnijgeneral}) is given by
\[
B_n\left( i,j\right) = \int _0^{c_i}k\left( \left( c_i-s\right) h_n\right) L_j\left( s\right) \, \mathrm{d}s.
\]
In this case, $B_n\left( i,j\right) $ is independent from $t_n$ and, given
some collocation parameters, it only depends on $h_n$.
\end{remark}

\begin{remark}
\label{nota0}
The coefficients $Z_{n,i}\geq 0$ given by (\ref{eqznigeneral}) are
positive, since $G$ is a positive function. But it does not imply that
$z_h$ was positive.
\end{remark}

The advantage of implicitly linear collocation methods (called
\textit{new collocation-type methods}, see \cite[p. 118]{Brun04}) lies
in the fact that, in contrast to direct collocation methods for
(\ref{eq0}), the integrals need not to be re-computed for every
iteration step when solving the nonlinear algebraic system
(\ref{eqznigeneral}).

\subsection{Nontrivial collocation solutions}

In this section we are going to recall the definition of
\textit{nontrivial solution} for the implicitly linear HVIE
(\ref{eq2}) with convolution kernel, and its corresponding collocation
problem. Nevertheless, it can be easily extended for the original
problem given by equation (\ref{eq0}),
and for nonconvolution kernels (see Remark
\ref{notanontrnonconv} below), but first we need the following
definition:

\begin{definition} We say that a property $\mathcal{P}$ holds
  \textit{near zero} if there exists $\epsilon >0$ such that
  $\mathcal{P}$ holds on $\left] 0,\delta \right[ $ for all $0<\delta
  <\epsilon $.  On the other hand, we say that $\mathcal{P}$ holds
  \textit{away from zero} if $\mathcal{P}$ holds on $\left] t,+\infty
  \right[ $ for all $t>0$.
\end{definition}

Given an implicitly linear HVIE (\ref{eq2}), the zero function is
always a solution, as it happens with equation (\ref{eq0}). Moreover,
for convolution kernels, given a solution $z(t)$ of (\ref{eq2}), and
$0<c<T$ the \textit{$c$-translated function} of $z$ given by
\[
z_c\left( t\right) :=\left\{
  \begin{array}{lll}
    0 & \, \,\mathrm{if} & 0\leq t<c \\
    \\
    z\left( t-c\right) & \, \, \mathrm{if} & c\leq t\leq T \\
\end{array}
\right.
\]
is also a solution of (\ref{eq2}). Thus, for convolution kernels, we
say that a solution is \textit{nontrivial} if it is neither
identically zero nor a $c$-translated function of another solution.

In this case,  $z$ is nontrivial if and only if it is
not identically zero near zero. This characterization allows us to
extend the concept of nontrivial solution to collocation problems with
convolution kernels:

\begin{definition}
\label{def:solcolnotri}
Given a collocation problem with convolution kernel, we say that a
collocation solution is \textit{nontrivial} if it is not identically
zero in $\sigma _0$.
\end{definition}

\begin{remark}
\label{notanontrnonconv}
The concept of nontrivial collocation solution can be easily extended
to nonconvolution kernels. Nevertheless, we have to take into account
that the $c$-translation of a solution of an implicitly linear HVIE
with nonconvolution kernel is not necessarily a solution, and there
can exist solutions that are $c$-translations of functions which are
not solutions.
\end{remark}

\section{Existence and uniqueness of nontrivial collocation solutions}
\label{sec:exist}

Given a kernel $K$, a nonlinearity $G$ and some collocation parameters
$\left\{ c_1,\ldots ,c_m\right\} $, our aim is to study the existence
of nontrivial collocation solutions (of the corresponding collocation
problem) in an interval $I=\left[ 0,T\right] $ using a mesh $I_h$. We
are only interested in existence (and uniqueness) properties where $h$
can be arbitrarily small and $N$ arbitrarily large, because
collocation solutions should converge to solutions of (\ref{eq2}) (if
they exist) when $h\rightarrow 0^+$ (and $N\rightarrow +\infty
$). Unfortunately, since uniqueness is not guaranteed, these
convergence problems are, in general, very complex and they are not in
the scope of this work. Taking this into account, we are going to
define three different kinds of existence of nontrivial collocation
solutions.

\begin{itemize}
\item We say that there is \textbf{existence near zero} if there
  exists $H_0>0$ such that if $0<h_0\leq H_0$ then there are
  nontrivial collocation solutions in $\left[ 0,t_1\right] $;
  moreover, there exists $H_n>0$ such that if $0<h_n\leq H_n$
  then there are nontrivial collocation solutions in $\left[
    0,t_{n+1}\right] $ (for $n=1,\ldots ,N-1$ and given $h_0,\ldots ,h_{n-1}>0$ such
  that there are nontrivial collocation solutions in $\left[
    0,t_n\right] $). Note that, in general, $H_n$ depends on
  $h_0,\ldots ,h_{n-1}$.

  This is the most general case of existence that we are going to
  consider. Loosely speaking this definition means that \emph{``collocation
  solutions can always be extended a bit more''}. Nevertheless, the
  existence of nontrivial collocation solutions for arbitrarily large
  $T$ is not assured (for instance, if there is a
  \textit{blow-up}). So, this is an existence \textit{near zero}.

\item We say that there is \textbf{existence for fine meshes} if there
  exists $H>0$ such that if $0<h\leq H$ then the corresponding
  collocation problem has nontrivial collocation solutions.

  This is a particular case of existence near zero, but in this case
  it is assured the existence of nontrivial collocation solutions in
  any interval $I$, using fine enough meshes $I_h$.

\item We say that there is \textbf{unconditional existence} if there
  exist nontrivial collocation solutions in any interval $I$ and for
  any mesh $I_h$.
\end{itemize}

We are going to study two cases of collocation problems:
\begin{itemize}
\item Case 1: $m=1$ with $c_1>0$.
\item Case 2: $m=2$ with $c_1=0$.
\end{itemize}
In these cases, the system (\ref{eqznigeneral}) is reduced to a single
nonlinear equation, whose solution is given by the fixed points of
$G\left( \alpha +\beta y\right) $ for some $\alpha, \beta $.

Let us state some lemmas which will be needed. Taking into account the
\textit{general conditions} over $G$, these results can be easily
proved (see Appendix A).

\begin{lemma}
  \label{lema1} The following statements are equivalent to the
  statement that $\frac{G\left( y\right) }{y}$ is unbounded (in
  $\left] 0,+\infty \right[ $):
\begin{itemize}
\item[(i)] There exists $\beta _0>0$ such that $G\left( \beta y\right)
  $ has nonzero fixed points for all $0<\beta \leq \beta _0$.
\item[(ii)] Given $A\geq 0$, there exists $\beta _A>0$ such that
  $G\left( \alpha +\beta y\right) $ has nonzero fixed points for all
  $0\leq \alpha \leq A$ and for all $0<\beta \leq \beta _A$.
\end{itemize}
\end{lemma}

\begin{lemma}
\label{lema2}
If $\frac{G\left( y\right) }{y}$ is unbounded near zero but it is
bounded away from zero, then there exists $\beta _0>0$ such that
$G\left( \alpha +\beta y\right) $ has nonzero fixed points for all
$\alpha \geq 0$ and for all $0<\beta \leq \beta _0$.

If, in addition, $\frac{G\left( y\right) }{y}$ is a strictly
decreasing function, then there exists a unique nonzero fixed point.
\end{lemma}

\begin{lemma}
\label{lema3}
If $\frac{G\left( y\right) }{y}$ is unbounded near zero and there
exists a sequence $\left\{ y_n\right\} _{n=1}^{+\infty }$ of positive
real numbers and divergent to $+\infty $ such that $\lim
_{n\rightarrow +\infty} \frac{G\left( y_n\right) }{y_n}=0$, then
$G\left( \alpha +\beta y\right) $ has nonzero fixed points for all
$\alpha \geq 0$ and for all $\beta >0$.

If, in addition, $\frac{G\left( y\right) }{y}$ is a strictly
decreasing function, then there exists a unique nonzero fixed point.
\end{lemma}

\begin{remark}
\label{notanonulo}
The nonzero fixed points are always strictly positive, since $G$ is
strictly positive in $\left] 0,+\infty \right[ $.  If, in addition,
$G\left( \alpha +\beta y\right) $ has fixed points for $\alpha >0$,
then these fixed points are necessarily nonzero and, hence, strictly
positive.
\end{remark}


So, we are going to study different kinds of existence and
uniqueness of nontrivial collocation solutions by means of the
equations (\ref{eqznigeneral}) considering two special
cases.

\subsection{Case 1: $m=1$ with $c_1>0$}

First, we shall consider $m=1$ with $c_1>0$. Note that if $m=1$ with
$c_1=0$, then the unique collocation solution is the trivial one.

We have the equations
\begin{equation}
\label{eqzn11}
Z_{n,1}=G\left( F_n\left( t_{n,1}\right) +h_n B_n Z_{n,1}\right)
\qquad (n=0,\ldots ,N-1),
\end{equation}
where
\begin{equation}
\label{eqbnnoconv}
B_n := B_n\left( 1,1\right) =\int _0^{c_1}K\left( t_{n,1},t_n+sh_n\right)
\, \mathrm{d}s,
\end{equation}
and the lag terms $F_n\left( t_{n,1}\right) $ are given by
(\ref{eqlag0}) with $i=1$.  From the \textit{general conditions}
imposed over $K$ it is assured that $B_n>0$.

\begin{remark}
\label{notahnbn}
Since $K$ is locally bounded, we have $h_nB_n\rightarrow 0$ when
$h_n\rightarrow 0^+$.
\end{remark}

Now, we are in position to give a characterization of the existence
near zero of nontrivial collocation solutions.

\begin{proposition}
\label{prop1}
Let $K$ be a kernel such that $K\left( t,s\right) \leq K\left( t',s\right) $ for all $0\leq s\leq
t<t'$.
Then there is existence near zero
if and only if $\frac{G\left( y\right) }{y}$ is unbounded.
\end{proposition}

Next, we are going to give some sufficient conditions for the
existence and uniqueness of nontrivial collocation solutions.

\begin{proposition}
\label{prop2}
If $\frac{G\left( y\right) }{y}$ is unbounded near zero but it is
bounded away from zero, then there is existence near zero.

If, in addition, $\frac{G\left( y\right) }{y}$ is a strictly
decreasing function, then there is at most one nontrivial collocation
solution.

For convolution kernels $K(t,s)=k(t-s)$, we can assure existence for fine meshes.
\end{proposition}

\begin{proposition}
\label{prop3}
If $\frac{G\left( y\right) }{y}$ is unbounded near zero and there
exists a sequence $\left\{ y_n\right\} _{n=1}^{+\infty }$ of positive
real numbers and divergent to $+\infty $ such that $\lim
_{n\rightarrow +\infty} \frac{G\left( y_n\right) }{y_n}=0$, then there
is unconditional existence.

If, in addition, $\frac{G\left( y\right) }{y}$ is a strictly
decreasing function, then there is at most one nontrivial collocation
solution.
\end{proposition}

In the proofs of Propositions \ref{prop1}, \ref{prop2} and \ref{prop3}
(see Appendix A) it is shown that the nontrivial collocation
solutions $z_h$ are strictly positive. Moreover, if $K\left( t,s\right) \leq K\left( t',s\right) $ for all $0\leq s\leq
t<t'$ and taking into account (\ref{eqcolloc}), it can
be easily proved that these collocation solutions are strictly
increasing functions.

\subsection{Case 2: $m=2$ with $c_1=0$}

Now, we are going to consider $m=2$ with $c_1=0$. Hence, we have to
solve the following equations:
\begin{eqnarray}
\label{eqzn12}
Z_{n,1}&=&G\left( F_n\left( t_{n,1}\right) \right) \\
\label{eqzn2}
Z_{n,2}&=&G\left( F_n\left( t_{n,2}\right) +h_nB_n\left( 2,1\right)
  Z_{n,1}+h_nB_n\left( 2,2\right) Z_{n,2}\right)
\quad (n=0,\ldots ,N-1),
\end{eqnarray}
where
\begin{equation}
\label{eqbnnoconv2}
B_n\left( 2,j\right) = \int _0^{c_2}K\left( t_{n,2},t_n+sh_n\right) L_j\left( s\right) \, \mathrm{d}s,\qquad (j=1,2)
\end{equation}
and $F_n\left( t_{n,i}\right) $ is given by (\ref{eqlag0}) with $m=2$.
Note that for solving $Z_{n,2}$ this system of equations can be
reduced to the single equation (\ref{eqzn2}), since (\ref{eqzn12})
gives us $Z_{n,1}$ directly.

From the \textit{general conditions} imposed over $K$ and taking into
account that functions $L_1\left( s\right) =1-\frac{s}{c_2}$ and
$L_2\left( s\right) =\frac{s}{c_2}$ are strictly positive in $\left]
  0,c_2\right[ $, it is assured that $B_n\left( 2,j\right) >0$ for
$j=1,2$.

\begin{remark}
\label{notahnbn2}
As in Remark \ref{notahnbn}, $h_nB_n\left( 2,j\right) \rightarrow 0$
when $h_n\rightarrow 0^+$. In particular, $h_nB_n\left( 2,2\right) \rightarrow 0$.
\end{remark}

Analogously to the previous case, we present similar results for
existence and uniqueness of nontrivial collocation solutions.

\begin{proposition}
\label{prop4}
Let $K$ be a kernel such that $K\left( t,s\right) \leq K\left( t',s\right) $ for all $0\leq s\leq
t<t'$. Then there is existence near zero
if and only if $\frac{G\left( y\right) }{y}$ is unbounded.
\end{proposition}

\begin{proposition}
\label{prop5}
Let $K$ be a kernel such that $K\left( t,s\right) \leq K\left( t',s\right) $ for all $0\leq s\leq
t<t'$.
If $\frac{G\left( y\right) }{y}$ is
unbounded near zero but it is bounded away from zero, then there is
existence near zero.

If, in addition, $\frac{G\left( y\right) }{y}$ is a strictly
decreasing function, then there is at most one nontrivial collocation
solution.

For convolution kernels $K(t,s)=k(t-s)$, the hypothesis on $K$ means that $k$ is increasing, and we can assure existence for fine meshes.
\end{proposition}

\begin{proposition}
\label{prop6}
Let $K$ be a kernel such that $K\left( t,s\right) \leq K\left( t',s\right) $ for all $0\leq s\leq
t<t'$. If $\frac{G\left( y\right) }{y}$ is
unbounded near zero and there exists a sequence $\left\{ y_n\right\}
_{n=1}^{+\infty }$ of positive real numbers and divergent to $+\infty
$ such that $\lim _{n\rightarrow +\infty} \frac{G\left( y_n\right)
}{y_n}=0$, then there is unconditional existence.

If, in addition, $\frac{G\left( y\right) }{y}$ is a strictly
decreasing function, then there is at most one nontrivial collocation
solution.
\end{proposition}

As in the previous case, the nontrivial collocation
solutions $z_h$ are strictly positive (see Appendix A). Also, from the
hypothesis on $K$ and (\ref{eqcolloc}), it can be easily
proved that the images under $z_h$ of the collocation points form a
strictly increasing sequence,
i.e. $0=Z_{0,1}<Z_{0,2}<Z_{1,1}<Z_{1,2}<\ldots
<Z_{N-1,2}$. Particularly, if $c_2=1$ then we can assure that $z_h$ is
a strictly increasing function.

\begin{remark}
\label{nota6}
In the case $c_2=1$, we can remove the hypothesis ``$K\left( t,s\right) \leq K\left( t',s\right) $ for all $0\leq s\leq
t<t'$'' from Propositions \ref{prop5} and \ref{prop6}, because we do not need
that $Z_{l,2}>Z_{l,1}$ with $l=1,\ldots ,n-1$ for proving that $z_h$
is strictly positive in $\left] 0,t_n\right[ $, and hence, we do not
need to prove that $Z_{n,2}>Z_{n,1}$ (see Appendix A). On the other hand, if we
remove this hypothesis, we can not assure that $z_h$ is a strictly
increasing function, as in the case $m=1$ studied in
Propositions \ref{prop2} and \ref{prop3}.
\end{remark}

\subsection{Nondivergent existence and uniqueness}
\label{sec:nodiv}

Given a kernel $K$, a nonlinearity $G$ and some collocation
parameters, we are interested in the study of existence of nontrivial
collocation solutions using meshes $I_h$ with arbitrarily small
$h>0$. Following this criterion, we are not interested in collocation
problems whose collocation solutions ``escape'' to $+\infty $ in a
certain $\sigma _n$ when $h_n\rightarrow 0^+$, since this is a
divergence symptom.

Let $S$ be an index set of all the nontrivial collocation solutions of
a collocation problem with mesh $I_h$. For any $s\in S$ we denote by
$Z_{s;n,i}$ the coefficients verifying equations (\ref{eqznigeneral})
(with $n=0,\ldots ,N-1$ and $i=1,\ldots ,m$) and such that, at least,
one of the coefficients $Z_{s;0,i}$ is different from zero (for some
$i\in\left\{ 1,\ldots ,m\right\} $).

So, given $K$, $G$ and some collocation parameters, we are going to
define the concepts of \textit{nondivergent existence} and
\textit{nondivergent uniqueness} of nontrivial collocation solutions.

\begin{definition}
  Given $0=t_0<\ldots <t_n$ such that there exist nontrivial
  collocation solutions using the mesh $t_0<\ldots <t_n$, we say that
  there is \textit{nondivergent existence in $t_n^+$} if
\[
\mathcal{Z}_{h_n}:=\inf _{s\in S_{h_n}} \left\{ \max _{i=1,\ldots
    ,m}\left\{ Z_{s;n,i}\right\} \right\}
\]
exists for small enough $h_n>0$ and it does not diverge to $+\infty $
when $h_n\rightarrow 0^+$.
\end{definition}

Note that the index set of the nontrivial collocation solutions is
denoted by $S_{h_n}$ because if we change $h_n$ then we change the
collocation problem.

\begin{definition}
  Given $I_h= \left\{ 0=t_0<\ldots<t_{N-1}\right\} $ such that there
  exist nontrivial collocation solutions using this mesh, we say that
  there is \textit{nondivergent existence} if there is nondivergent
  existence in $t_n^+$ for $n=0,\ldots, N-1$.
\end{definition}

For studying the concept of \textit{nondivergent uniqueness} we need to
state the following definitions.

\begin{definition}
  Given $0=t_0<\ldots <t_n$ such that there exist nontrivial
  collocation solutions using the mesh $t_0<\ldots <t_n$, we say that
  there is \textit{nondivergent uniqueness in $t_n^+$} if
\[
\min _{s\in S_{h_n}} \left\{ \max _{i=1,\ldots
    ,m}\left\{ Z_{s;n,i}\right\} \right\}=\mathcal{Z}_{h_n}
\]
exists for small enough $h_n>0$, and it does not diverge to $+\infty $
when $h_n\rightarrow 0^+$, but
\[
\inf _{s\in S_{h_n}} \left( \left\{ \max _{i=1,\ldots ,m}\left\{
      Z_{s;n,i}\right\} \right\} -\left\{ \mathcal{Z}_{h_n}\right\}
\right)
\]
diverges (note that $\inf \emptyset = +\infty $).
\end{definition}

\begin{definition}
  Given $I_h= \left\{ 0=t_0<\ldots<t_{N-1}\right\} $ such that there
  exist nontrivial collocation solutions using this mesh, we say that
  there is \textit{nondivergent uniqueness} if there is nondivergent
  uniqueness in $t_n^+$ for $n=0,\ldots, N-1$.
\end{definition}

When ``nondivergent uniqueness'' is assured, but there is not
``uniqueness'' (of nontrivial collocation solutions), there is only
one nontrivial collocation solution that makes sense, as it is stated
in the following definition.

\begin{definition}
  In case of nondivergent uniqueness, the \textit{nondivergent
    collocation solution} is the one whose coefficients $Z_{n,i}$
  satisfy $\max _{i=1,\ldots ,m}\left\{ Z_{n,i}\right\}
  =\mathcal{Z}_{h_n}$ for $n=0,\ldots ,N-1$.
\end{definition}

Next, we are going to study nondivergent existence and uniqueness
for cases 1 ($m=1$ with $c_1>0$) and 2 ($m=2$ with $c_1=0$). Recall
that the \textit{general conditions} over $K$ and $G$ are always held,
even if it is not explicitly mentioned.

But first, we are going to state a result that reduces the study of
nondivergent existence (and uniqueness) for any mesh $I_h$ to the
study of nondivergent existence (and uniqueness) in $t_0^+$. Lemma
\ref{lemma29} can be easily proved taking into account the
\textit{general conditions} over $G$.

\begin{lemma}
\label{lemma29}
Given $\alpha >0$, the minimum of the nonzero fixed points of $G\left(
  \alpha +\beta y\right) $ exists for a small enough $\beta >0$ and
converges to $G\left( \alpha \right) $ when $\beta \rightarrow 0^+$.

If, in addition, $\frac{G\left( \alpha +y\right) }{y}$ is strictly
decreasing near zero, then the other nonzero fixed points (if they
exist) diverge to $+\infty $.
\end{lemma}

As a consequence of this lemma, we obtain the following proposition.

\begin{proposition}
\label{prop:1nodiv}
In cases 1 and 2 with existence of nontrivial collocation solutions,
``nondivergent existence'' and ``nondivergent existence in $t_0^+$''
are equivalent.

If, in addition, $\frac{G\left( \alpha +y\right) }{y}$ is strictly
decreasing near zero for all $\alpha >0$, then ``nondivergent
uniqueness'' and ``nondivergent uniqueness in $t_0^+$'' are
equivalent.
\end{proposition}

\begin{remark}
\label{nota31}
The condition ``$\frac{G\left( \alpha +y\right) }{y}$ strictly
decreasing near zero for all $\alpha >0$'' is very weak, since $\lim
_{y\rightarrow 0^+} \frac{G\left( \alpha +y\right) }{y}=+\infty $. For
example, if $G$ is twice differentiable a.e.  without accumulation of
non-differentiable points and without accumulation of sign changes of
the second derivative, then this condition is held. So, if a
collocation problem has nontrivial collocation solutions and, roughly
speaking, $G$ is ``well-behaved", then nondivergent uniqueness in
$t_0^+$ implies nondivergent uniqueness in any $t_n^+$.
\end{remark}

Next, we are going to give a characterization of nondivergent existence, but
first we need the following lemma, that can be easily proved
taking into account the \textit{general conditions} over $G$.

\begin{lemma}
\label{lem:2nodiv}
If $\frac{G\left( y\right) }{y}$ is unbounded, then the minimum of the
nonzero fixed points of $G\left( \beta y\right) $ exists for a small
enough $\beta >0$. In this case, this minimum does not diverge to
$+\infty $ when $\beta \rightarrow 0^+$ if and only if $\frac{G\left(
y\right) }{y}$ is unbounded near zero.

If, in addition, $\frac{G\left( y\right) }{y}$ is strictly decreasing
near zero, then the other nonzero fixed points (if they exist) diverge
to $+\infty $.
\end{lemma}

Using Lemma \ref{lem:2nodiv} and Proposition \ref{prop:1nodiv}, it can
be proved the next result.

\begin{proposition}
\label{prop7}
In cases 1 and 2 with existence of nontrivial collocation solutions,
there is nondivergent existence if and only if $\frac{G\left( y\right)
}{y}$ is unbounded near zero.

If, in addition, $\frac{G\left( \alpha +y\right) }{y}$ is strictly
decreasing near zero for all $\alpha \geq 0$, then there is
nondivergent uniqueness.
\end{proposition}

\begin{remark}
\label{remweird2}
If $\frac{G\left( y\right) }{y}$ is unbounded near zero, then the
condition ``$\frac{G\left( \alpha +y\right) }{y}$ strictly decreasing
near zero for all $\alpha \geq 0$'' is very weak, using the same
arguments as in Remark \ref{nota31}. So, if $G$ is ``well-behaved''
(see Remark \ref{nota31}), nondivergent existence implies also
nondivergent uniqueness.
\end{remark}

To sum up, combining Proposition \ref{prop7}
with Propositions \ref{prop1}, \ref{prop2}, \ref{prop3} (case 1) and
\ref{prop4}, \ref{prop5}, \ref{prop6} (case 2) and taking into account
Remark \ref{nota6}:

\begin{itemize}

\item $K\left( t,s\right) \leq K\left( t',s\right) $ for all $0\leq s\leq t<t'$;

$\frac{G\left( y\right) }{y}$ is unbounded
  near zero $\Leftrightarrow $ Nondivergent existence near zero.

  Moreover, if $G$ is ``well-behaved''
  $\Rightarrow $ Nondivergent uniqueness near zero.

\item (Hypothesis only for case 2: $c_2=1$, or $K\left( t,s\right) \leq K\left( t',s\right) $ for all $0\leq s\leq t<t'$);

  $\frac{G\left( y\right) }{y}$ is unbounded near zero but bounded
  away from zero $\Rightarrow $ Nondivergent existence near zero.

  Moreover, if $G$ is ``well-behaved''
  $\Rightarrow $ Nondivergent uniqueness near zero.

  For convolution kernels $K(t,s)=k(t-s)$ $\Rightarrow $ Nondivergent existence or uniqueness (resp.) for fine meshes.

\item (Hypothesis only for case 2: $c_2=1$, or $K\left( t,s\right) \leq K\left( t',s\right) $ for all $0\leq s\leq t<t'$);

  $\frac{G\left( y\right) }{y}$ is unbounded near zero and there
  exists a sequence $\left\{ y_n\right\} _{n=1}^{+\infty }$ of
  positive real numbers and divergent to $+\infty $ such that $\lim
  _{n\rightarrow +\infty} \frac{G\left( y_n\right) }{y_n}=0$
  $\Rightarrow $ Unconditional nondivergent existence.

  Moreover, if $G$ is ``well-behaved''
  $\Rightarrow $ Unconditional nondivergent uniqueness.

\end{itemize}

\section{Examples}
\label{sec:examples}

\subsection{Numerical study of convergence}
\label{sec:numconvergence}

In this section we are going to show an example of a collocation
problem (with nondivergent uniqueness in cases 1 and 2), and study
numerically how the nontrivial collocation solution $z_h$ converges to
a solution of the implicitly linear HVIE (\ref{eq2}) when
$h\rightarrow 0^+$. In fact, we are going to show how the function
$y_h=\mathcal{V}z_h$ (see Remark \ref{remnlHVIE}) converges to a
solution of the original nonlinear HVHIE (\ref{eq0}).

We are going to consider the equation
\begin{equation}
\label{eqex1}
y\left( t\right) = \int _0 ^t \left( t-s\right) \left( y\left( s\right) \right)^{1/2} \,
\mathrm{d}s , \qquad t\in \left[ 0,1\right] .
\end{equation}
The kernel and the nonlinearity verify the \textit{general conditions}
and all the hypotheses of Propositions \ref{prop3} (for case 1) and
\ref{prop6} (for case 2). Hence there is unconditional existence and
uniqueness; moreover, we can also apply Proposition \ref{prop7} for assuring that the unique nontrivial collocation solution is
nondivergent.

\begin{figure}[tbp]
\begin{center}
  \includegraphics[width=1\textwidth]{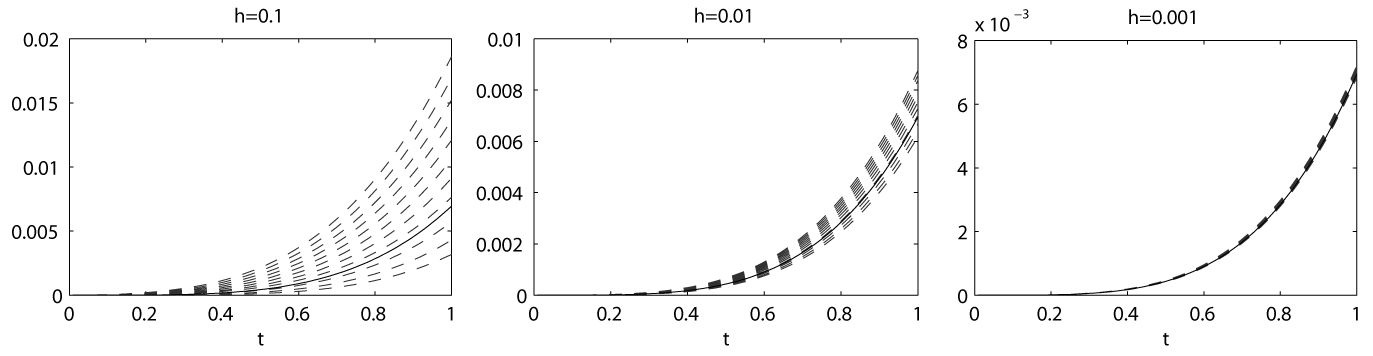}
\end{center}
\caption{Case 1 ($m=1$, $c_1>0$). Approximation $y_h$ (grey dots) for
  different stepsizes $h$, with $c_1$ from $0.01$ (lower) to $1$
  (upper). The solution $y$ is also represented (black solid line).}
\label{m1_soluciones}
\end{figure}

\begin{figure}[tbp]
\begin{center}
  \includegraphics[width=1\textwidth]{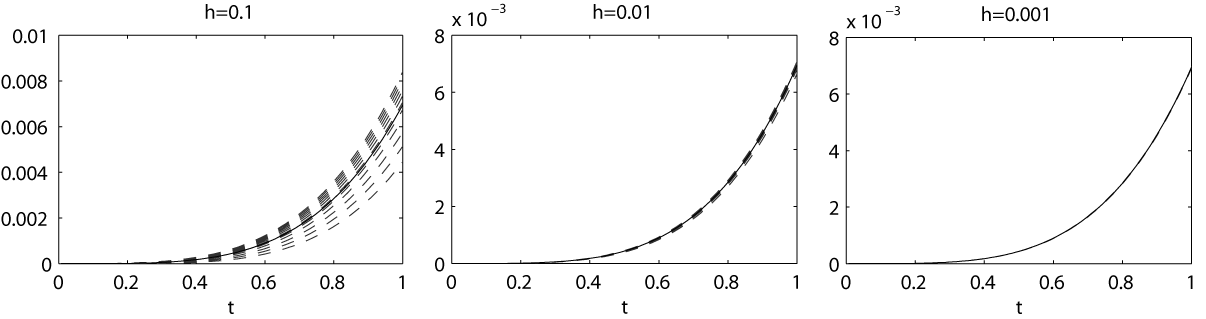}
\end{center}
\caption{Case 2 ($m=2$, $c_1=0$). Approximation $y_h$ (grey dots) for
  different stepsizes $h$, with $c_2$ from $0.01$ (lower) to $1$
  (upper). The solution $y$ is also represented (black solid line).}
\label{m2_soluciones}
\end{figure}

Since the unique nontrivial solution of the nonlinear HVHIE (\ref{eqex1}) is given
by
\[
y\left( t\right) =\frac{1}{144}t^4,
\]
we study the difference between this solution and $y_h$ when
$h\rightarrow 0^+$ (see Figures \ref{m1_soluciones} and
\ref{m2_soluciones}).

\begin{figure}[tbp]
\begin{center}
  \includegraphics[width=1\textwidth]{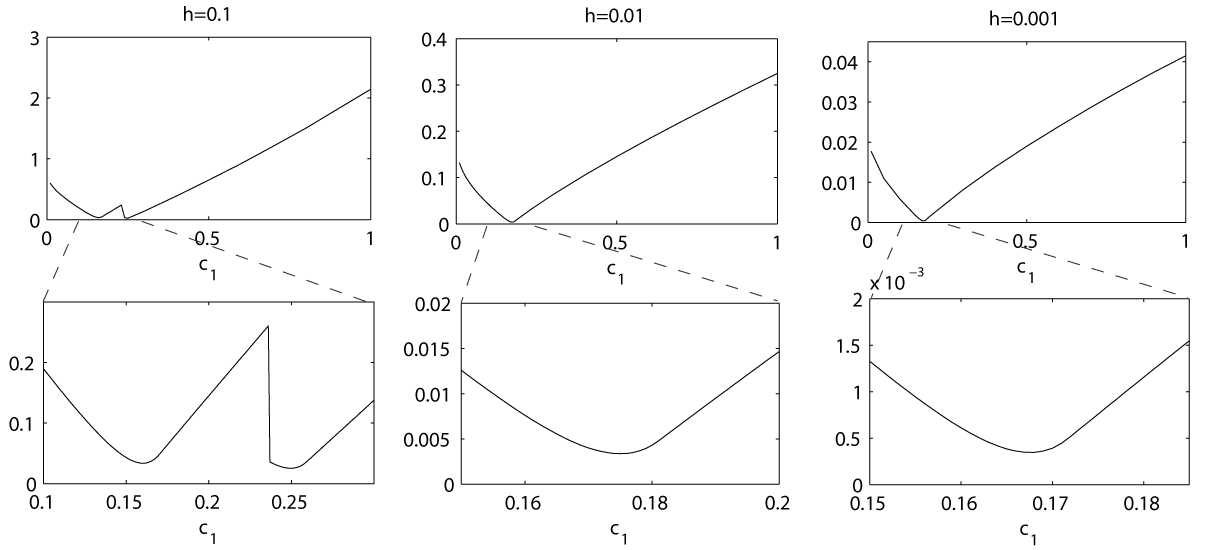}
\end{center}
\caption{Case 1. Relative error varying $c_1$ (from $0.01$ to $1$) for
  different stepsizes $h$.}
\label{m1_yh-y}
\end{figure}

\begin{figure}[tbp]
\begin{center}
  \includegraphics[width=1\textwidth]{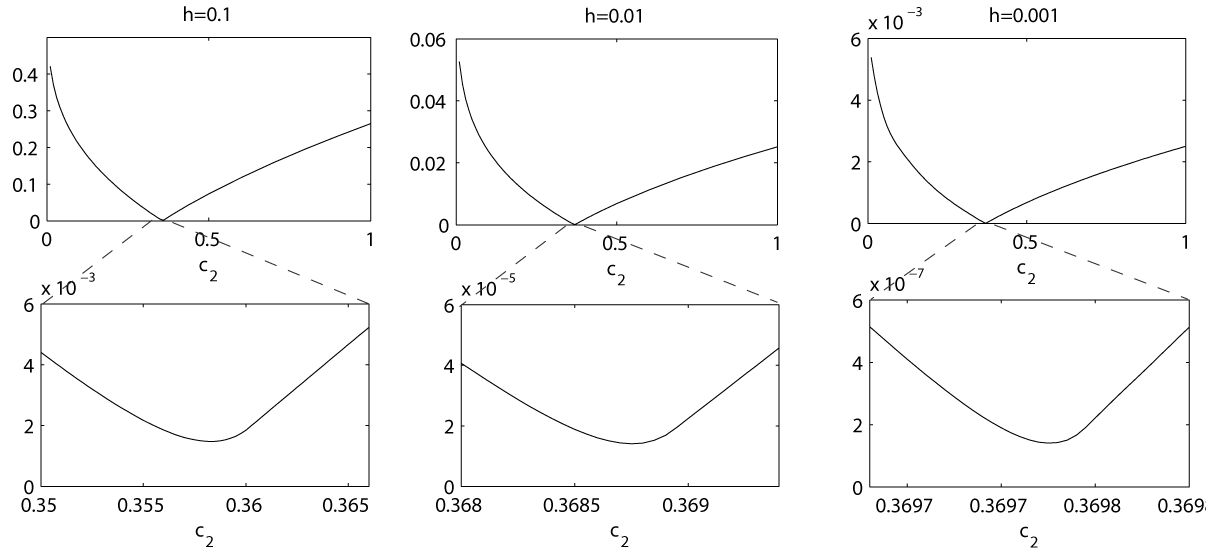}
\end{center}
\caption{Case 2. Relative error varying $c_2$ (from $0.01$ to $1$) for
  different stepsizes $h$.}
\label{m2_yh-y}
\end{figure}

In Figures \ref{m1_yh-y} and \ref{m2_yh-y} it is represented the
relative error $\frac{\int _I \left| y_h-y\right| }{\int _I y}$
varying $c_1$ (case 1) or $c_2$ (case 2), for different stepsizes
$h$. It is shown that the rate of convergence is the same as
$h$. Moreover, a convenient choice of the collocation parameter can
reduce the relative error more than two orders, as we see in Tables
\ref{table1} and \ref{table2}. Nevertheless, how to find a good collocation parameter is an open problem.

\subsection{Examples in other cases}

\begin{table}
\begin{center}
\begin{tabular}{|l|c|c|c|}
  \hline
  Case 1 ($m=1$, $c_1>0$) & $h=0.1$ & $h=0.01$ & $h=0.001$ \\
  \hline
  Max. rel. error & $2.1$ & $3.2\cdot 10^{-1}$ & $4.1\cdot 10^{-2}$ \\
  & \small{($c_1=1$)} & \small{($c_1=1$)} & \small{($c_1=1$)} \\
  \hline
  Min. rel. error & $2.5\cdot 10^{-2}$ & $3.4\cdot 10^{-3}$ & $3.4\cdot 10^{-4}$ \\
  & \small{($c_1=0.25$)} & \small{($c_1=0.175$)} & \small{($c_1=0.168$)} \\
  \hline
\end{tabular}
\end{center}
\caption{Case 1. Maximum and minimum relative errors varying $c_1$, for different stepsizes $h$ (see Figure \ref{m1_yh-y}).}
\label{table1}
\end{table}

\begin{table}
\begin{center}
\begin{tabular}{|l|c|c|c|}
  \hline
  Case 2 ($m=2$, $c_1=0$) & $h=0.1$ & $h=0.01$ & $h=0.001$ \\
  \hline
  Max. rel. error & $4.2\cdot 10^{-1}$ & $5.2\cdot 10^{-2}$ & $5.4\cdot 10^{-3}$ \\
  & \small{($c_2=0.01$)} & \small{($c_2=0.01$)} & \small{($c_2=0.01$)} \\
  \hline
  Min. rel. error & $1.5\cdot 10^{-3}$ & $1.4\cdot 10^{-5}$ & $1.4\cdot 10^{-7}$ \\
  & \small{($c_2=0.358$)} & \small{($c_2=0.369$)} & \small{($c_2=0.37$)} \\
  \hline
\end{tabular}
\end{center}
\caption{Case 2. Maximum and minimum relative errors varying $c_2$, for different stepsizes $h$ (see Figure \ref{m2_yh-y}).}
\label{table2}
\end{table}

If the collocation problem is not in the scope of cases 1 and 2, the
existence of nontrivial collocation solutions is not assured, even if
the kernel and the nonlinearity satisfy all the mentioned conditions.

For example, if we consider
\begin{itemize}
\item $K\left( t,s\right) =\left( t-s\right) ^a,\quad a>0$,
\item $G\left( y\right) = y^{1/b}, \quad b>1$,
\end{itemize}
the \textit{general conditions} over $K$ and $G$ are held, as well as
all the conditions for assuring unconditional existence and uniqueness
in cases 1 and 2. Moreover, by Proposition \ref{prop7}, the unique
nontrivial collocation solution is nondivergent. Nevertheless, in
other cases we can not assure anything, as we will show in the next
examples (where existence is considered at least in $\sigma _0$).

\begin{itemize}

\item Let $b=2$. If $m=2$ and $c_1>0$ then it can be proved that there
  is not any nontrivial collocation solution (for any $a>0$).

  But, if $m=3$, $a=2$ and $c_1>0$ then there exists a unique
  nontrivial collocation solution at least in $\sigma _0$,
  independently of $c_2$ and $c_3$; on the other hand, if $c_1=0$ then
  there is not any nontrivial collocation solution (contrary to what
  happens in the case $m=2$).

  Moreover, if $m=3$ and $a=1$ then there exists a unique nontrivial
  collocation solution at least in $\sigma _0$, independently of the
  collocation parameters.

\item Let $b=3$, $a=1$ and $m\geq 2$. It can be proved that there
  exists a unique nontrivial collocation solution at least in $\sigma
  _0$, independently of the collocation parameters. This is a special
  case, since the collocation solution for $m\geq 2$ coincides with
  the solution of the corresponding implicitly linear HVIE, $z\left(
    t\right) =t/\sqrt{6}$.
\end{itemize}

\section{Discussion and comments}

A general theoretical analysis of numerical approximations of
nontrivial collocation solutions for equation (\ref{eq0}) is an
outstanding problem which, to our knowledge, remains open. Among the
reasons for such difficulty is the lack of uniqueness of solutions
and the lack of Lipschitz-continuity of the nonlinearity, to name a
few.

Nevertheless, this work is a first step towards such general analysis:
we give conditions for different kinds of existence and uniqueness of
nontrivial collocation solutions aiming also for the convergence
analysis (e.g. nondivergent existence).

We have also studied numerically a concrete example showing the
accuracy of the methods and how their errors depend on the collocation
points (see Section \ref{sec:numconvergence}). Therefore, in the case
of nondivergent existence and uniqueness of collocation solutions,
collocation methods have proved to be a valuable numerical tool for
approximating solutions. Moreover, these methods are strengthened by
the fact that, if $\frac{G\left( y\right) }{y}$ is a strictly
decreasing function, the nonzero fixed points mentioned in Lemma
\ref{lema2} and Lemma \ref{lema3} are attractors, so they can be
easily found numerically with iteration techniques. Furthermore, if
there is nondivergent existence, the minimum of such nonzero fixed
points is also an attractor.

It is worth pointing out that the conditions we have imposed over the
nonlinearity in order to assure the nondivergent existence near zero
of nontrivial collocation solutions, namely $\frac{G\left( y\right)
}{y}$ must be unbounded near zero, agrees with the conditions for
existence of nontrivial solutions of the original equation
(\ref{eq0}). To be more precise, if $\frac{G\left( y\right) }{y}$ is
unbounded near zero then it is not Lipschitz-continuous, which is a
necessary condition for the existence of nontrivial solutions of
equation (\ref{eq0}).

\appendix
\section{Proofs}
\label{sec:app}

\begin{pol1}

($\Rightarrow (i)$) Let us prove that if $\frac{G\left( y\right) }{y}$ is
unbounded (in $\left] 0,+\infty \right[ $), then there exists $\beta _0>0$ such that $G\left( \beta y\right) $
has nonzero fixed points for all $0<\beta \leq \beta _0$.
Let us take $y_0>0$, $y_1:=G\left( y_0\right) >0$ and $\beta _0:=\frac{y_0}{y_1}>0$. So,
given $0<\beta <\beta _0$, we have
\begin{equation}
\label{plema1.1}
G\left( \beta y_1\right) <G\left( \beta _0 y_1\right) =G\left( y_0 \right) =y_1,
\end{equation}
since $G$ is strictly increasing. Moreover, by hypothesis, there exists $y_2>0$ such that $\frac{G\left( y_2\right) }{y_2}>\frac{1}{\beta }$. Let us define $y_3:=\frac{y_2}{\beta }$, then
\begin{equation}
\label{plema1.2}
G\left( \beta y_3\right) = G\left( y_2\right) >\frac{y_2}{\beta }=y_3.
\end{equation}
So, taking into account (\ref{plema1.1}) and (\ref{plema1.2}), $G\left( \beta y\right) $ has fixed points between $y_1$ and $y_3$, since it is continuous.

($\Leftarrow (i)$) Let us prove the other implication. Given $M>\frac{1}{\beta _0}>0$, we take $0<\beta <\frac{1}{M}<\beta _0$. Then, by hypothesis, there exists $y_0>0$ such that $G\left( \beta y_0\right) =y_0$, and so, taking $y_1:=\beta y_0$ we have $\frac{G\left( y_1\right) }{y_1}=\frac{1}{\beta}>M$.

($(i)\Rightarrow (ii)$) Let us prove that $(i)$ implies $(ii)$, for $A>0$. Let us take $y_0>0$, $y_1:=G\left( A+y_0\right) >0$ and $\beta _A:=\frac{y_0}{y_1}>0$. So,
given $0<\beta <\beta _A$ and $0<\alpha \leq A$ (we can suppose $\alpha >0$ taking $\beta _A\leq \beta _0$), we have
\begin{equation}
\label{plema1.3}
G\left( \alpha+\beta y_1\right) <G\left( A+\beta _A y_1\right) =G\left( A+y_0 \right) =y_1,
\end{equation}
since $G$ is strictly increasing. Moreover
\begin{equation}
\label{plema1.4}
G\left( \alpha +\beta 0\right) = G\left( \alpha \right) >0.
\end{equation}
So, taking into account (\ref{plema1.3}) and (\ref{plema1.4}), $G\left( \alpha +\beta y\right) $ has nonzero fixed points between $0$ and $y_1$, since it is continuous.

($(ii)\Rightarrow (i)$) Trivial.
\end{pol1}

\begin{pol2}
Given $\alpha \geq 0$ and a small enough $\beta >0$, it is easy to prove (taking analogous arguments as in the proof of Lemma \ref{lema1}) that there exists $y_1>0$ such that $G\left( \alpha +\beta y_1\right) <y_1$, since $\frac{G\left( y\right) }{y}$ is bounded away from zero. On the other hand, it is easy to prove that there exists $y_3\geq 0$ such that $G\left( \alpha +\beta y_3\right) >y_3$, since $\frac{G\left( y\right) }{y}$ is unbounded near zero (note that the case $\alpha >0$ is trivial taking $y_3=0$). So, $G\left( \alpha +\beta y\right) $ has nonzero fixed points between $y_3$ and $y_1$, since it is continuous.
Moreover, it is clear that if $\frac{G\left( y\right) }{y}$ is a strictly decreasing function, this fixed point is unique.
\end{pol2}

\begin{pol3}
Analogous to the proof of Lemma \ref{lema2}, but taking any $\beta >0$.
\end{pol3}

\begin{pot1}

($\Leftarrow $) Let us prove that if $\frac{G\left( y\right) }{y}$
is unbounded, then there is existence near zero. So, we are going to
prove by induction over $n$ that there exist $H_n>0$ $\left(
n=0,\ldots ,N-1\right) $ such that if $0<h_n\leq H_n$ then there
exist solutions of the system (\ref{eqzn11}) with $Z_{0,1}>0$:

\begin{itemize}

\item For $n=0$, taking into account Remark \ref{notahnbn} and Lemma
  \ref{lema1}-\textit{(i)}, we choose a small enough $H_0>0$ such that
  $0<h_0B_0\leq \beta _0$ for all $0<h_0\leq H_0$. So, since the lag
  term is $0$, we can apply Lemma \ref{lema1}-\textit{(i)} to the
  equation (\ref{eqzn11}), concluding that there exist strictly
  positive solutions for $Z_{0,1}$.

\item Let us suppose that, choosing one of those $Z_{0,1}$, there
  exist $H_1,\ldots ,H_{n-1}>0$ such that if $0<h_i\leq H_i$
  ($i=1,\ldots ,n-1$) then there exist coefficients $Z_{1,1},\ldots
  ,Z_{n-1,1}$ fulfilling the equation (\ref{eqzn11}). Note that these
  coefficients are strictly positive by Remark \ref{notanonulo}, and
  hence, it is guaranteed that the corresponding collocation solution
  $z_h$ is strictly positive in $\left] 0,t_n\right[ $.

\item Finally, we are going to prove that there exists $H_n>0$ such
  that if $0<h_n\leq H_n$ then there exists $Z_{n,1}>0$ fulfilling the
  equation (\ref{eqzn11}) with the previous coefficients
  $Z_{0,1},\ldots Z_{n-1,1}$:

  Let us define $A:={F_n\left( t_{n}+c_1\right) }$. So, taking into
  account Remark \ref{notahnbn} and Lemma \ref{lema1}-\textit{(ii)},
  we choose a small enough $0<H_n\leq 1$ such that $0<h_nB_n\leq \beta
  _A$ for all $0<h_n\leq H_n$.  Then, we have $0<F_n\left(
    t_{n,1}\right) \leq F_n\left( t_{n}+c_1\right) =A$ because the hypothesis over $K$ (and the \textit{general
    conditions}), $z_h$ is strictly positive in $\left] 0,t_n\right[
  $, and $t_{n,1}\leq t_n+c_1$. Hence, we can apply Lemma
  \ref{lema1}-\textit{(ii)}, obtaining the existence of $Z_{n,1}$
  (that is strictly positive by Remark \ref{notanonulo}).

\end{itemize}

($\Rightarrow $) For proving the other condition, we use Lemma
\ref{lema1}-\textit{(i)}, taking into account Remark \ref{notahnbn}.
\end{pot1}

\begin{pot2}
  First we are going to consider a convolution kernel $K(t,s)=k(t-s)$.
  Taking into account Remarks \ref{notabnconv}, \ref{notahnbn}, and
  Lemma \ref{lema2}, we choose a small enough $H>0$ such that
\begin{equation}
\label{eqhint0c1}
h\int _0^{c_1}k\left( \left( c_1-s\right) h\right) \, \mathrm{d}s\leq \beta _0
\end{equation}
for all $0<h\leq H$.

So, if $0<h\leq H$, we are going to prove by induction over $n$ that
there exist solutions of the system (\ref{eqzn11}) with $Z_{0,1}>0$:

\begin{itemize}

\item For $n=0$, by (\ref{eqhint0c1}) we have $0<h_0B_0\leq \beta _0$,
  because $0<h_0\leq h\leq H$. Since the lag term is $0$, we can apply
  Lemma \ref{lema2} to the equation (\ref{eqzn11}), concluding that
  there exist strictly positive solutions for $Z_{0,1}$.

\item Let us suppose that, choosing one of those $Z_{0,1}$, there
  exist coefficients $Z_{1,1},\ldots ,Z_{n-1,1}$ fulfilling the
  equation (\ref{eqzn11}). Note that these coefficients are strictly
  positive by Remark \ref{notanonulo}, and hence, it is guaranteed
  that the corresponding collocation solution $z_h$ is strictly
  positive in $\left] 0,t_n\right[ $.

\item Finally, we are going to prove that there exists $Z_{n,1}>0$
  fulfilling the equation (\ref{eqzn11}) with the previous
  coefficients $Z_{0,1},\ldots Z_{n-1,1}$:

  On one hand, taking into account (\ref{eqlag0}), the lag term
  $F_n\left( t_{n,1}\right) $ is strictly positive because $z_h$ is
  strictly positive in $\left] 0,t_n\right[ $ and $k$ satisfies the
  \textit{general conditions}. On the other hand, $0<h_nB_n\leq \beta
  _0$ because $0<h_n\leq h\leq H$. Hence, we can apply Lemma
  \ref{lema2} to the equation (\ref{eqzn11}), obtaining the existence
  of $Z_{n,1}$ (that is strictly positive by Remark \ref{notanonulo}).

\end{itemize}

If, in addition, $\frac{G\left( y\right) }{y}$ is a strictly
decreasing function, the uniqueness is an immediate consequence of
Lemma \ref{lema2}.

For general (nonconvolution) kernels, the proof is analogous but
  $B_n$ is given by (\ref{eqbnnoconv}) and it does not only depend on
  $h_n$. The existence of $h_n$ lies in choosing a small enough $H_n$
  such that $h_nB_n\leq \beta _0$ for all $0<h_n\leq H_n$.
\end{pot2}

\begin{pot3}
  The proof is analogous to the proof of Proposition \ref{prop2}, but
  in this case $H>0$ is any positive real number, and using Lemma
  \ref{lema3} instead of Lemma \ref{lema2}.
\end{pot3}

\begin{pot4}

  ($\Leftarrow $) Let us prove that if $\frac{G\left( y\right) }{y}$
  is unbounded, then there is existence near zero. So, we are going to
  prove by induction over $n$ that there exist $H_n>0$ $\left(
    n=0,\ldots ,N-1\right) $ such that if $0<h_n\leq H_n$ then there
  exist solutions of the system (\ref{eqzn2}) with $Z_{0,2}>0$:

\begin{itemize}

\item For $n=0$, taking into account Remark \ref{notahnbn2} and Lemma
  \ref{lema1}-\textit{(i)}, we choose a small enough $H_0>0$ such that
  $0<h_0B_0\left( 2,2\right) \leq \beta _0$ for all $0<h_0\leq
  H_0$. So, since the lag terms are $0$ and $Z_{0,1}=G\left( 0\right)
  =0$, we can apply Lemma \ref{lema1}-\textit{(i)} to the equation
  (\ref{eqzn2}), concluding that there exist strictly positive
  solutions for $Z_{0,2}$.

\item Let us suppose that, choosing one of those $Z_{0,2}$, there
  exist $H_1,\ldots ,H_{n-1}>0$ such that if $0<h_i\leq H_i$
  ($i=1,\ldots ,n-1$) then there exist coefficients $Z_{1,2},\ldots
  ,Z_{n-1,2}$ fulfilling the equation (\ref{eqzn2}). Moreover, let us
  suppose that these coefficients satisfy $Z_{l,2}>Z_{l,1}>0$ for
  $l=1,\ldots ,n-1$, and hence, it is guaranteed that the
  corresponding collocation solution $z_h$ is strictly positive in
  $\left] 0,t_n\right[ $.

\item Finally, we are going to prove that there exists $H_n>0$ such
  that if $0<h_n\leq H_n$ then there exists $Z_{n,2}>0$ fulfilling the
  equation (\ref{eqzn2}) with the previous coefficients, and
  $Z_{n,2}>Z_{n,1}>0$:

    Let us define
    \[
    A:=F_n\left( t_n+c_2\right) +\int _0^{c_2}k\left( \left(
        c_2-s\right) \right) L_1\left( s\right) \, \mathrm{d}s \,\,
    G\left( F_n\left( t_n+c_1\right) \right) .
    \]
    So, taking into account Remark \ref{notahnbn2} and Lemma
    \ref{lema1}-\textit{(ii)}, we choose a small enough $0<H_n\leq 1$
    such that $0<h_nB_n\left( 2,2\right) \leq \beta _A$ for all
    $0<h_n\leq H_n$.  Taking into account (\ref{eqlag0}), the lag
    terms $F_n\left( t_{n,i}\right) $ are strictly positive for
    $i=1,2$, because $z_h$ is strictly positive in $\left]
      0,t_n\right[ $ and $K$ satisfies the \textit{general
      conditions}.  Therefore, by (\ref{eqzn12}), $Z_{n,1}=G\left(
      F_n\left( t_{n,1}\right) \right) $ is strictly positive, because
    $G$ is strictly positive in $\left] 0,+\infty \right[ $.
    Moreover, $h_nB_n\left( 2,1\right) >0$, and hence, $h_nB_n\left(
      2,1\right) Z_{n,1}>0$.  So,
    \[
    0<F_n\left( t_{n,2}\right) +h_nB_n\left( 2,1\right) Z_{n,1}\leq A,
    \]
    because $K$ satisfies the hypothesis, $z_h$ is strictly positive
    in $\left] 0,t_n\right[ $, the polynomial $L_1\left( s\right) $ is
    strictly positive in $\left] 0,c_2\right[ $, the nonlinearity $G$
    is a strictly increasing function, and $t_{n,j}\leq t_n+c_j$ for
    $j=1,2$.  Hence, we can apply Lemma \ref{lema1}-\textit{(ii)} to
    the equation (\ref{eqzn2}), obtaining the existence of $Z_{n,2}$.

    Concluding, we have to check that $Z_{n,2}>Z_{n,1}$. Since $K$ satisfies the hypothesis in (\ref{eqlag0}), $F_n\left(
      t_{n,2}\right) \geq F_n\left( t_{n,1}\right) $, and hence, by
    the properties of $G$, we have
    \[
    \begin{array}{lll}
      Z_{n,2}&=&G\left( F_n\left( t_{n,2}\right) +h_n\sum_{j=1}^2B_n\left( 2,j\right) Z_{n,j}\right)\\
      \\
      &>&G\left( F_n\left( t_{n,2}\right) \right) \, \, \geq \, \, G\left( F_n\left( t_{n,1}\right)
      \right) \, \, =\, \, Z_{n,1}.
    \end{array}
    \]

\end{itemize}

($\Rightarrow $) For proving the other condition, we use Lemma
\ref{lema1}-\textit{(i)}, taking into account Remark \ref{notahnbn2}.
\end{pot4}

\begin{pot5}
  First we are going to consider a convolution kernel $K(t,s)=k(t-s)$.
  Taking into account Remarks \ref{notabnconv}, \ref{notahnbn2}, and
  Lemma \ref{lema2}, we choose a small enough $H>0$ such that
\begin{equation}
  h\int _0^{c_2}k\left( \left( c_2-s\right) h\right) \frac{s}{c_2}\, \mathrm{d}s\leq \beta _0
\end{equation}
for all $0<h\leq H$. From here, the proof is analogous to the proof
of Proposition \ref{prop4}, where $H_n=H$ and using Lemma \ref{lema2}
instead of Lemma \ref{lema1}. So, we do not need any $A$, and we do
not need to check that $F_n\left( t_{n,2}\right) +h_nB_n\left(
  2,1\right) Z_{n,1}\leq A$.

If, in addition, $\frac{G\left( y\right) }{y}$ is a strictly
decreasing function, the uniqueness is an immediate consequence of Lemma
\ref{lema2}.

For general (nonconvolution) kernels, the proof is analogous but
  $B_n\left( 2,j\right) $ is given by (\ref{eqbnnoconv2}) and it does
  not only depend on $h_n$. The existence of $h_n$ lies in choosing a
  small enough $H_n$ such that $h_nB_n\left( 2,2\right) \leq \beta _0$
  for all $0<h_n\leq H_n$.
\end{pot5}

\begin{pot6}
  The proof is analogous to the proof of Proposition \ref{prop5}, but
  in this case $H>0$ is any positive real number, and using Lemma
  \ref{lema3} instead of Lemma \ref{lema2}.
\end{pot6}

\newpage
\thispagestyle{empty}
\mbox{}

\end{document}